\documentclass[preprint,12pt]{elsarticle}
\usepackage{amssymb,amsmath,amsthm}

\biboptions{comma,square,compress}

\journal{Indagationes Mathematicae}

\renewcommand*{\leq}{\leqslant} \renewcommand*{\geq}{\geqslant}

\newcommand*{\vare}{\varepsilon} \newcommand*{\varp}{\varphi}

\newcommand*{\bw}{\mathbf w}

\newcommand*{\mR}{\mathbb R} \newcommand*{\mS}{\mathbb S}
\newcommand*{\mT}{\mathbb T} \newcommand*{\mZ}{\mathbb Z}

\newcommand*{\cG}{\mathcal G} \newcommand*{\cH}{\mathcal H}
\newcommand*{\cK}{\mathcal K} \newcommand*{\cL}{\mathcal L}
\newcommand*{\cM}{\mathcal M} \newcommand*{\cN}{\mathcal N}
\newcommand*{\cO}{\mathcal O} \newcommand*{\cP}{\mathcal P}
\newcommand*{\cQ}{\mathcal Q} \newcommand*{\cS}{\mathcal S}
\newcommand*{\cT}{\mathcal T} \newcommand*{\cV}{\mathcal V}
\newcommand*{\cY}{\mathcal Y}

\newcommand*{\fZ}{\mathfrak Z}

\newcommand*{\rmd}{\mathrm d} \newcommand*{\rmi}{\mathrm i}
\newcommand*{\const}{\mathrm{const}}

\newcommand*{\Ad}{\mathop{\mathrm{Ad}}\nolimits}
\newcommand*{\codim}{\mathop{\mathrm{codim}}\nolimits}
\newcommand*{\Fix}{\mathop{\mathrm{Fix}}\nolimits}
\newcommand*{\GL}{\mathop{\mathrm{GL}}\nolimits}
\newcommand*{\rank}{\mathop{\mathrm{rank}}\nolimits}

\newcommand*{\gl}{\mathop{\mathfrak{gl}}\nolimits}

\newcommand*{\Landau}[1]{\mathrm O\bigl(#1\bigr)}

\newcommand*{\skobki}[1]{\textup{(}#1\textup{)}}
\newcommand*{\kvadra}[1]{\textup{[}#1\textup{]}}

\theoremstyle{plain}
\newtheorem{thm}{Theorem}
\newtheorem*{Herlem}{Herman's lemma}
\newtheorem*{reflem}{Standard reflection lemma}
\newtheorem*{Hampar}{Hamiltonian KAM paradigm (for exact symplectic forms)}
\newtheorem*{revpar}{Reversible KAM paradigm}

\theoremstyle{definition}
\newtheorem{dfn}{Definition}
\newtheorem{rem}{Remark}

\begin{document}

\begin{frontmatter}

\title{Quasi-periodic perturbations within the reversible context~2 \\
in KAM theory}

\author{Mikhail B. Sevryuk}

\ead{sevryuk@mccme.ru}

\address{Institute of Energy Problems of Chemical Physics,
The Russia Academy of Sciences, \\
Leninski\u{\i} prospect~38, Bldg.~2, Moscow 119334, Russia}

\begin{abstract}
The paper consists of two sections. In Section~1, we give a short review of
KAM theory with an emphasis on Whitney smooth families of invariant tori in
typical Hamiltonian and reversible systems. In Section~2, we prove a KAM-type
result for non-autonomous reversible systems (depending quasi-periodically on
time) within the almost unexplored reversible context~2. This context refers
to the situation where $\dim\Fix G<\frac{1}{2}\codim\cT$, here $\Fix G$ is the
fixed point manifold of the reversing involution $G$ and $\cT$ is the
invariant torus one deals with.
\end{abstract}

\begin{keyword}
KAM theory \sep Reversible systems \sep Hamiltonian systems \sep
Whitney smooth Cantor-like families of quasi-periodic invariant tori \sep
Quasi-periodic perturbations \sep Reducibility \sep
Artificial additional parameter

\MSC[2010] 70K43 \sep 70H33 \sep 37J40 \sep 70H08
\end{keyword}

\end{frontmatter}

\section{KAM theory from a bird's eye view}

KAM (Kolmogorov--Arnold--Moser) theory is the theory of quasi-periodic motions
(i.e.\ conditionally periodic motions with incommensurable frequencies) in
non-integrable dynamical systems. The phase curves of such motions fill up
densely invariant tori (so called quasi-periodic invariant tori) in the phase
space. In turn, these tori are usually organized into complicated hierarchical
structures consisting of tori of different dimensions. However, the ``building
blocks'' of such structures are Whitney smooth Cantor-like families of
invariant tori rather than individual tori. We refer the reader to
e.g.\ \cite{BHS96Gro,BHS96LNM} for a precise definition of Whitney smooth
families of invariant tori (to be brief, a Whitney smooth function on a
closed set in $\mR^\ell$ is a function extendible to a smooth function on a
neighborhood of this set). The properties of Whitney smooth families of
invariant tori depend strongly on the symmetries preserved by the system in
question.

\subsection{Hamiltonian systems}

As our first example, consider autonomous Hamiltonian vector fields on a
finite dimensional manifold equipped with an \emph{exact} symplectic $2$-form.
In the Hamiltonian KAM theory, the following easy observation is of crucial
importance.

\begin{Herlem}
Any quasi-periodic invariant torus of a Hamiltonian system is isotropic
provided that the symplectic form is exact. In particular, the dimension of
such a torus does not exceed the number of degrees of freedom.
\end{Herlem}

Recall that a submanifold $\cL$ of a symplectic manifold is said to be
\emph{isotropic} if the restriction of the symplectic form to $\cL$ vanishes,
i.e., the tangent space $T_a\cL$ is contained in its skew orthogonal
complement for each point $a\in\cL$. Herman himself \cite{H88,H89} proved the
lemma above in a certain particular case but the general case
\cite{BHS96LNM,S03} is not harder at all (in fact, this lemma goes back to
Moser \cite[pp.~157--158]{M66}). Herman's lemma can be carried over to locally
Hamiltonian systems, i.e.\ to Hamiltonian systems with multi-valued Hamilton
functions (for instance, the system $\dot{x}=0$, $\dot{y}=1$ on a cylinder
with coordinates $x\in\mS^1$, $y\in\mR$ and the symplectic form
$\rmd x\wedge\rmd y$ is locally Hamiltonian: it is afforded by the Hamilton
function $x$ with values in $\mS^1$ rather than in $\mR$).

\begin{rem}\label{dva}
Of course, invariant tori of dimensions $0$ (equilibria) and $1$ (periodic
trajectories) are always isotropic whether or not the symplectic form is
exact. The same is valid for invariant $2$-tori. Indeed, consider an invariant
manifold $\cL$ of a Hamiltonian flow with Hamilton function $\cH$. Suppose
that $\cH|_{\cL}=\const$ and that almost all the points of $\cL$ are not
equilibria. It is not hard to verify that $\cL$ is isotropic if $\dim\cL=2$
and coisotropic if $\codim\cL=2$ (regardless of whether the symplectic form is
exact).
\end{rem}

Recall that a submanifold $\cL$ of a symplectic manifold is said to be
\emph{coisotropic} if the tangent space $T_a\cL$ contains its skew orthogonal
complement for each point $a\in\cL$.

The following concepts are also of principal importance in KAM theory.

\begin{dfn}
Let an invariant $n$-torus $\cT$ of some flow on an $(n+\ell)$-dimensional
manifold carry conditionally periodic motions with frequency vector
$\omega\in\mR^n$. This torus is said to be \emph{reducible} (or
\emph{Floquet}) if in a neighborhood of $\cT$, there exists a coordinate frame
$x\in\mT^n$, $w\in\cO_\ell(0)$ in which the torus $\cT$ itself is given by the
equation $\{w=0\}$ and the dynamical system takes the \emph{Floquet form}
$\dot{x}=\omega+\Landau{|w|}$, $\dot{w}=Lw+\Landau{|w|^2}$ with an
$x$-independent matrix $L\in\gl(\ell,\mR)$. This matrix (not determined
uniquely) is called the \emph{Floquet matrix} of the torus $\cT$, and its
eigenvalues are called the \emph{Floquet exponents} of $\cT$.
\end{dfn}

Here and henceforth, $\mT^n=(\mS^1)^n=(\mR/2\pi\mZ)^n$ is the standard
$n$-torus, while $\cO_\ell(a)$ denotes an unspecified neighborhood of a point
$a\in\mR^\ell$.

In other words, an invariant torus is reducible if the variational equations
along this torus can be reduced to a form with constant coefficients.

The essence of the Hamiltonian KAM theory can now be formulated as follows.

\begin{Hampar}
In a \emph{typical} Hamiltonian system with $M\geq 1$ degrees of freedom
\skobki{in the case of an exact symplectic $2$-form}, there are
\begin{itemize}
\item isolated equilibria,
\item smooth one-parameter families of closed trajectories \skobki{one
trajectory per energy value},
\item and Whitney smooth Cantor-like $n$-parameter families of isotropic
invariant $n$-tori carrying quasi-periodic motions with strongly
incommensurable \skobki{e.g.\ Diophantine} frequencies for each
$2\leq n\leq M$.
\end{itemize}
These tori can be either reducible or non-reducible. The Floquet exponents of
a reducible invariant $n$-torus include value $0$ of multiplicity $n$
\skobki{for $1\leq n\leq M$} whereas the remaining $2M-2n$ Floquet exponents
come in pairs $\lambda,-\lambda$ \skobki{for $0\leq n\leq M-1$}.
\end{Hampar}

The word ``typical'' here means that Hamiltonian systems with the properties
indicated constitute an open set (to be more precise, a set with non-empty
interior) in the functional space of all the Hamiltonian systems with $M$
degrees of freedom. The meaning of this word in the sequel will be similar.

Isotropic invariant $n$-tori in a Hamiltonian system with $M$ degrees of
freedom are said to be \emph{Lagrangian} for $n=M$ and \emph{lower
dimensional} for $n<M$.

By now, Whitney smooth Cantor-like families of isotropic invariant tori in
Hamiltonian systems have been thoroughly explored, especially in the reducible
case. The reader is referred to
\cite{AKN06,BHS96Gro,BHS96LNM,BHT90,BS10,QS93,S95RMS,S03} for surveys, precise
statements, and bibliographies. Here we will confine ourselves with the
following remark. In KAM theory, there are known phenomena leading to a
decrease or increase in the dimension of invariant tori.

One of the ``lowering'' phenomena is destruction of resonant tori. Consider a
partially integrable Hamiltonian system possessing a smooth $n$-parameter
family of isotropic invariant $n$-tori ($n\geq 2$) carrying conditionally
periodic motions with frequency vectors $\omega(\mu)$, where $\mu$ is the
parameter of the family. Let $1\leq r\leq n-1$. Typically, this $n$-parameter
family of $n$-tori contains an $(n-r)$-parameter subfamily of tori whose
frequencies satisfy $r$ independent \emph{fixed} resonance relations
$\left\langle q^{(i)},\omega(\mu)\right\rangle=0$, $q^{(i)}\in\mZ^n$,
$1\leq i\leq r$ (here and henceforth, the angle brackets denote the standard
inner product). Then, under a generic Hamiltonian perturbation, this smooth
$(n-r)$-parameter subfamily of resonant $n$-tori gives rise to a finite
collection of Whitney smooth Cantor-like $(n-r)$-parameter families of
isotropic quasi-periodic invariant $(n-r)$-tori (so called Treshch\"ev tori,
see \cite{AKN06,BHS96LNM,BS10,S03} for surveys and references). Of course, for
$r=n-1$, these $1$-parameter families of closed trajectories (called
Poincar\'e trajectories) are smooth rather than Cantor-like. In fact, break-up
of resonant tori has been studied by now in the Lagrangian case only (for $n$
equal to the number $M$ of degrees of freedom).

One of the ``raising'' phenomena is excitation of elliptic normal modes.
Consider an $n$-parameter family of reducible isotropic invariant $n$-tori in
a Hamiltonian system with $M>n$ degrees of freedom. Let $\omega(\mu)$ be the
frequency vectors of the tori, where $\mu$ is the parameter of the family.
Suppose that among the non-zero Floquet exponents of each of these tori,
there are $r$ pairs ($1\leq r\leq M-n$) of purely imaginary numbers
$\pm\rmi\beta_1(\mu),\ldots,\pm\rmi\beta_r(\mu)$ that depend smoothly on $\mu$
(in the Whitney sense for $n\geq 2$). The remaining $2(M-n-r)$ non-zero
Floquet exponents are also allowed to be purely imaginary. Then, generically,
in a neighborhood of this $n$-parameter family of $n$-tori, there is an
$(n+r)$-parameter family of reducible isotropic invariant $(n+r)$-tori. The
frequencies of these tori are close to
$\omega_1(\mu),\ldots,\omega_n(\mu),\beta_1(\mu),\ldots,\beta_r(\mu)$.
One says that the Floquet exponents
$\pm\rmi\beta_1(\mu),\ldots,\pm\rmi\beta_r(\mu)$ ``excite'', see again
\cite{AKN06,BHS96LNM,BS10,S03} for surveys and references. For instance, a
neighborhood of a generic elliptic equilibrium ($n=0$, $r=M$) contains a
Whitney smooth Cantor-like $M$-parameter family of Lagrangian invariant
$M$-tori. The simplest example is $n=0$, $r=M=1$: an elliptic equilibrium of
a planar Hamiltonian system is surrounded by a smooth one-parameter family of
periodic trajectories (energy levels). In the case where $n=0$, $r=1$, $M>1$,
similar smooth one-parameter families of closed trajectories are called
Lyapunov families.

We emphasize that, starting from $n$-parameter families of invariant $n$-tori,
we obtain $(n\pm r)$-parameter families of invariant $(n\pm r)$-tori, in a
complete agreement with the paradigm above.

In the case of non-exact symplectic forms, there become possible coisotropic
invariant $n$-tori with $M+1\leq n\leq 2M-1$ (such tori are studied in
Parasyuk's theory and are sometimes said to be \emph{higher dimensional}) as
well as so called \emph{atropic} invariant $n$-tori (i.e., tori that are
neither isotropic nor coisotropic) with $3\leq n\leq 2M-3$. As before, $M$ is
the number of degrees of freedom. Coisotropic and a~fortiori atropic invariant
tori have been explored much worse than isotropic ones, we refer the reader to
\cite{BHS96Gro,BHS96LNM,BS10,S03,S08} for discussions and relevant
bibliographies. Atropic quasi-periodic invariant tori of dimensions $2$ and
$2M-2$ cannot exist according to Remark~\ref{dva}.

\begin{rem}
In \cite{QS93}, we introduced the concept of $s$-exact symplectic forms: a
symplectic $2$-form is said to be \emph{$s$-exact} ($s\geq 1$) if its $s$th
exterior power is exact (it is clear that an $s$-exact symplectic form is also
$s'$-exact for any $s'>s$). Seizing the opportunity, we now make an important
comment to this concept. Consider a $2M$-dimensional connected manifold $\cM$.
If $\cM$ is \emph{closed} (i.e.\ compact without boundary) then it cannot
admit even $M$-exact symplectic forms since the $M$th exterior power of a
symplectic form on $\cM$ is a volume element. On the other hand, if $\cM$ is
\emph{open} (i.e.\ either noncompact or with a non-empty boundary) and carries
a nondegenerate $2$-form then it also admits an exact symplectic form. In
fact, according to Gromov's theorem (see e.g.\ \cite[Section~7.3]{MS98}),
every homotopy class of nondegenerate $2$-forms on $\cM$ can be represented by
a symplectic form representing any prescribed cohomology class in
$H^2(\cM,\mR)$. In particular, any manifold $\mR^k\times\mT^l$ with $k+l$ even
and $k\geq 1$ admits an exact symplectic form (this is of course obvious for
$k\geq l$ but not obvious at all for $0<k<l$).
\end{rem}

\subsection{Reversible systems}\label{revers}

Our second and main example is reversible vector fields. Consider a finite
dimensional connected manifold $\cM$ and a smooth involution $G:\cM\to\cM$
($G^2$ is the identity transformation).

\begin{dfn}\label{obrat}
A vector field $V$ on $\cM$ is said to be \emph{reversible} with respect to
involution $G$ (or $G$-reversible) if $\Ad G(V)=TG(V\circ G)=-V$ or,
equivalently, if the function $t\mapsto G\bigl(a(-t)\bigr)$ is a solution of
the equation $\dot{a}=V(a)$, $a\in\cM$, whenever $\mR\ni t\mapsto a(t)$ is.
\end{dfn}

For instance, the Newtonian equations of motion $\ddot{\bw}=F(\bw,\dot{\bw})$,
$\bw\in\mR^\ell$, are reversible with respect to the phase space involution
$G:(\bw,\dot{\bw})\mapsto(\bw,-\dot{\bw})$ if and only if the forces $F$ are
even in the velocities $\dot{\bw}$ (e.g.\ are independent of $\dot{\bw}$). The
papers \cite{LR98,RQ92} present general surveys of the theory of reversible
systems with extensive bibliographies.

If a submanifold $\cL\subset\cM$ is invariant under a $G$-reversible flow, so
is its ``mirror'' image $G(\cL)$. If $G(\cL)\neq\cL$ then, as a rule, the
dynamics near each of the two invariant manifolds $\cL$ and $G(\cL)$
considered separately exhibits no special features. Therefore, while speaking
of an invariant manifold $\cL$ of a $G$-reversible system, one usually assumes
$\cL$ to be invariant not only under the flow itself but also under the
reverser $G$. In the sequel, the words ``an invariant torus (in particular, an
equilibrium or periodic trajectory) of a reversible system'' will always mean
a torus invariant under both the flow and the reversing involution.

The crucial role in the dynamics of $G$-reversible systems on $\cM$ is played
by the fixed point set $\Fix G=\bigl\{a\in\cM \bigm| G(a)=a\bigr\}$. This set
is a submanifold of the same smoothness class as the involution $G$ itself
(a very particular case of Bochner's theorem \cite{B72,MZ74}, see also
\cite{CF64}). The fixed point manifold $\Fix G$ can well be empty or consist
of several connected components of different dimensions. Extensive information
on the fixed point submanifolds of involutions of various manifolds is
presented in e.g.\ the books \cite{B72,CF64}, see also the articles
\cite{QS93,S95RMS,S11}.

As is widely known, there is a deep similarity between reversible and
Hamiltonian dynamics \cite{A84,AS86,LR98,RQ92,S86,S91}. In particular, many
fundamental results of the Hamiltonian KAM theory possess reversible
counterparts. The reversible analog of Herman's lemma is the following (also
very easy) statement.

\begin{reflem}[\cite{BHS96Gro,BHS96LNM,S86,Snew}]
In any quasi-periodic invariant $n$-torus $\cT$ of a $G$-reversible flow, one
can introduce a coordinate frame $\varp\in\mT^n$ in which the dynamics on
$\cT$ takes the form $\dot{\varp}=\omega$ \skobki{$\omega\in\mR^n$ being the
frequency vector of $\cT$} and the restriction of $G$ to $\cT$ takes the form
$G|_{\cT}:\varp\mapsto-\varp$. In particular, $\cT\cap\Fix G$ consists of
$2^n$ points $\varp$ with $\varp_i$ equal to either $0$ or $\pi$,
$1\leq i\leq n$, and the dimension of any connected component of $\Fix G$
having a non-empty intersection with $\cT$ does not exceed $\codim\cT$.
Moreover, in a smooth \skobki{or Whitney smooth} family of quasi-periodic
invariant tori, the coordinate $\varp$ can be chosen to depend smoothly
\skobki{respectively Whitney smoothly} on the torus.
\end{reflem}

This lemma can be carried over to \emph{weakly reversible} systems,
i.e.\ systems reversed by phase space diffeomorphisms $G$ that are not
necessarily involutions \cite{A84,AS86,S86,Snew}. The definition of weakly
$G$-reversible vector fields $V$ has the form $\Ad G(V)=TG(V\circ G^{-1})=-V$.

In the sequel, we will consider only involutions $G:\cM\to\cM$ for which
$\Fix G\neq\varnothing$ and all the connected components of $\Fix G$ are of
the same dimension, so that $\dim\Fix G$ is well defined (this is the case for
almost all the reversible systems encountered in practice). We will write
$\dim\Fix G=P$ and $\codim\Fix G=Q$, so that $\dim\cM=P+Q$. According to the
standard reflection lemma, if a $G$-reversible system admits a quasi-periodic
invariant $n$-torus then $n\leq Q$.

To formulate the reversible counterpart of the Hamiltonian KAM paradigm, one
has to deal with smooth \emph{families of reversible systems} (depending on
an external $s$-dimensional parameter $\nu$) rather than with individual
reversible systems (corresponding to the case where $s=0$). The reversible
analog of the Hamiltonian KAM paradigm is the following statement.

\begin{revpar}
In the product of the phase space and the parameter space of a \emph{typical}
$s$-parameter family of $G$-reversible systems with $\dim\Fix G=P$ and
$\codim\Fix G=Q$, there are
\begin{itemize}
\item smooth $(P-Q+s)$-parameter families of equilibria
\skobki{for $s\geq\max\{Q-P, \, 0\}$},
\item smooth $(P-Q+s+1)$-parameter families of closed trajectories
\skobki{for $s\geq\max\{Q-P-1, \, 0\}$},
\item and Whitney smooth Cantor-like $(P-Q+s+n)$-parameter families of
invariant $n$-tori carrying quasi-periodic motions with strongly
incommensurable \skobki{e.g.\ Diophantine} frequencies for each
$2\leq n\leq Q$ \skobki{for $s\geq\max\{Q-P-n+1, \, 0\}$}.
\end{itemize}
These tori can be either reducible or non-reducible. The Floquet exponents of
a reducible invariant $n$-torus include value $0$ of multiplicity $|P-Q+n|$
\skobki{for $P\neq Q-n$} whereas the remaining $2\min\{P, \, Q-n\}$ Floquet
exponents come in pairs $\lambda,-\lambda$
\skobki{for $\min\{P, \, Q-n\}\geq 1$}.
\end{revpar}

In the case where $P>0$ and $Q>n$, the invariant $n$-tori in question are said
to be \emph{lower dimensional}.

Let us explain $|P-Q+n|$ zero eigenvalues of the Floquet matrix here. Consider
a matrix $L$ anti-commuting with a fixed involutive matrix
$\cK\in\GL(\cP+\cQ,\mR)$ [$\cK^2$ is the identity matrix], the eigenvalues $1$
and $-1$ of $\cK$ being of multiplicities $\cP$ and $\cQ$, respectively. If
$\cP\neq\cQ$ then $0$ is an eigenvalue of $L$ of multiplicity at least
$|\cP-\cQ|$ \cite{S86,S92}. The remaining $\cP+\cQ-|\cP-\cQ|=2\min\{\cP,\cQ\}$
eigenvalues of $L$ come in pairs $\lambda,-\lambda$ and are generically other
than zero \cite{S86,S92}.

\begin{rem}
One cannot hope to encounter an isolated quasi-periodic invariant torus of
dimension $n\geq 2$ in a generic dissipative, volume preserving, Hamiltonian,
or reversible system (or in a generic family of systems). Of course, generic
systems may admit isolated invariant $n$-tori with $n\geq 2$, but such tori
would not carry conditionally periodic motions (the induced dynamics would be
``phase-locked''). Therefore, for $(P-Q+s+n)$-parameter families of invariant
tori of dimensions $n\geq 2$ in the reversible KAM paradigm, one has
$P-Q+s+n\geq 1$, i.e.\ $s\geq Q-P-n+1$ rather than $s\geq Q-P-n$.
\end{rem}

Similarly to the Hamiltonian context, one may consider destruction of resonant
tori. Starting with a smooth $(P-Q+s+n)$-parameter family of invariant
$n$-tori carrying conditionally periodic motions, one would obtain Whitney
smooth Cantor-like $(P-Q+s+n-r)$-parameter families of quasi-periodic
invariant $(n-r)$-tori (for $1\leq r\leq n-1$ and $P-Q+s+n-r\geq 1$; it
suffices to require $P-Q+s+n-r=P-Q+s+1\geq 0$ for $r=n-1$). Excitation of
elliptic normal modes makes sense as well. Starting with a
$(P-Q+s+n)$-parameter family of reducible invariant $n$-tori, one expects
to obtain $(P-Q+s+n+r)$-parameter families of invariant $(n+r)$-tori for
$1\leq r\leq\min\{P, \, Q-n\}$.

However, while the Hamiltonian KAM paradigm has been proven completely by now,
one cannot say the same about the reversible KAM paradigm. Consider a
quasi-periodic invariant $n$-torus $\cT$ of a $G$-reversible system with
$\dim\Fix G=P$ and $\codim\Fix G=Q$. For this torus, we have two non-negative
``characteristic numbers'' $P$ and $Q-n$. The codimension of $\cT$ is their
sum $P+Q-n$.

\begin{dfn}
The situation where
\[
P\geq Q-n \quad\Longleftrightarrow\quad \dim\Fix G\geq\tfrac{1}{2}\codim\cT
\]
is called \emph{the reversible context~1}. The opposite situation where
\[
P<Q-n \quad\Longleftrightarrow\quad \dim\Fix G<\tfrac{1}{2}\codim\cT
\]
is called \emph{the reversible context~2}.
\end{dfn}

It turns out that by now, almost all the reversible KAM theory has been
devoted exclusively to the reversible context~1. This context has been nearly
as developed as the Hamiltonian KAM theory. Surveys, precise statements of the
theorems, and bibliographies are given in
e.g.\ \cite{AS86,BH95,BHS96Gro,BHS96LNM,M67,QS93,S86,S91,S95RMS,S95Cha,S98},
the reader is also referred to
\cite{BCH07,BCHV09,BHN07,L01,S06,S07DCDS,S07Stek,WX09,WXZ10,W01} for some
important results obtained after the review \cite{S98}. By the way, the
reversible context~1 requires an external parameter only in the case where
$n\geq 2$ and $P=Q-n$ (so that $s\geq 1$). For other values of $P$, $Q$, and
$n$ within the reversible context~1, one may consider individual systems
($s=0$).

The reversible context~2 was first described in \cite{BHS96Gro,BHS96LNM} where
the paradigm for this context was formulated as a conjecture. The paradigm for
the reversible context~1 was given in \cite{BHS96Gro,BHS96LNM} separately (as
we saw above, the paradigms for both the contexts can be unified). The task of
developing the reversible KAM theory in context~2 was listed (as problem~9)
among the ten problems of the classical KAM theory in the note \cite{S08}. The
first result in the reversible context~2 was obtained no earlier than in 2011
\cite{S11}. It concerned the ``\emph{extreme}'' reversible context~2 where
$P=0$. Lower dimensional reducible invariant tori in the reversible context~2
(with $0<P<Q-n$) were treated in \cite{Snew}. Both the papers \cite{S11,Snew}
examine only analytic families of invariant tori (in the analytic category,
i.e., under the assumption that the reversing involution, the vector fields
themselves and their families are analytic) in the presence of many external
parameters. To the best of the author's knowledge, all the KAM theory for the
reversible context~2 is currently confined to these two papers (and the
present one).

Such a situation seems in fact somewhat strange because the reversible
contexts~1 and~2 are closely related. Destruction of resonant tori allows one
to pass from context~1 to context~2 \cite{S11}: it is quite possible that
$P\geq Q-n$ but $P<Q-(n-r)$. Nevertheless, break-up of resonant tori in
reversible systems (for $r\leq n-2$) has been studied by now only in the case
where $P=Q=n$ \cite{L01,W01} and the inequality $P<Q-n+r$ is therefore never
met. Excitation of elliptic normal modes allows one to pass from context~2 to
context~1: it is quite possible that $P<Q-n$ but $P\geq Q-(n+r)$.

\begin{rem}
The statements ``matrices $L$ and $\cK$ anti-commute'' and ``matrices $L$ and
$-\cK$ anti-commute'' are equivalent. Thus, on the level of linear operators,
there is no difference between the reversible contexts~1 and~2.
\end{rem}

All the discussion above has been devoted to autonomous flows (either
Hamiltonian or reversible). The papers \cite{S11,Snew} constituting the first
steps in the reversible context~2 did not handle non-autonomous systems
either. In the next section, we will examine reversible systems (within
context~2) depending quasi-periodically on time. In context~1, such reversible
systems were dealt with in \cite{BCHV09,M65,M66,S07DCDS}. By the way, Moser's
note \cite{M65} is the first paper on the reversible KAM theory whatsoever.
Our paper \cite{S07DCDS} treats quasi-periodic perturbations in the reversible
context~1, the Hamiltonian context, the volume preserving context, and the
dissipative context from a unified viewpoint.

At the end of \cite{S11}, we listed ten tentative topics and directions for
further research. The eighth topic was non-autonomous perturbations depending
on time periodically or quasi-periodically.

Note that families of reducible quasi-periodic invariant tori in KAM theory
are Cantor-like because of resonances among the frequencies as well as of
those between the frequencies and the imaginary parts of the Floquet
exponents. Along an analytic subfamily of such tori, the frequencies and the
imaginary parts of the Floquet exponents should therefore be constants (at
least up to proportionality).

\section{Quasi-periodic perturbations within the reversible context~2}

For time-dependent vector fields on a manifold $\cM$ equipped with an
involution $G:\cM\to\cM$, Definition~\ref{obrat} of reversibility is modified
as follows.

\begin{dfn}
A time-dependent vector field $V_t$ on $\cM$ is said to be \emph{reversible}
with respect to involution $G$ (or $G$-reversible) if
$TG(V_t\circ G)\equiv-V_{-t}$ or, equivalently, if the function
$t\mapsto G\bigl(a(-t)\bigr)$ is a solution of the equation $\dot{a}=V_t(a)$,
$a\in\cM$, whenever $t\mapsto a(t)$ is.
\end{dfn}

On $\cM$, we will consider reversible systems $\dot{a}=V_t(a)$ depending on
time quasi-periodically with $N\geq 1$ incommensurable basic frequencies
$\Omega_1,\ldots,\Omega_N$:
\[
V_t(a)\equiv\cV(a,\Omega_1t,\ldots,\Omega_Nt),
\]
where the function $\cV=\cV(a,X_1,\ldots,X_N)$ is $2\pi$-periodic in each of
the arguments $X_1,\ldots,X_N$. From the viewpoint of KAM theory, the natural
problem is to look for quasi-periodic invariant tori of the corresponding
autonomous system
\begin{equation}
\dot{a}=\cV(a,X_1,\ldots,X_N), \quad \dot{X}=\Omega
\label{eqauto}
\end{equation}
on $\cM\times\mT^N$. Among the frequencies of such tori, there are $N$ numbers
$\Omega_1,\ldots,\Omega_N$. It is easy to see that a non-autonomous system
$\dot{a}=V_t(a)$ on $\cM$ is reversible with respect to involution $G$ if and
only if the autonomous system~\eqref{eqauto} on $\cM\times\mT^N$ is reversible
with respect to the involution $\cG:(a,X)\mapsto\bigl(G(a),-X\bigr)$.

Now we are in the position to state the main result of this paper. We will
consider analytic families of reducible quasi-periodic invariant tori in
quasi-periodic perturbations of autonomous reversible systems. Let
$x\in\mT^n$, $y\in\cY\subset\mR^m$, and $z\in\cO_{2p}(0)$ be the phase space
variables where $n\geq 0$, $m\geq 1$, $p\geq 0$, and $\cY$ is an open domain.
The reversing involution is $G:(x,y,z)\mapsto(-x,-y,Kz)$ where
$K\in\GL(2p,\mR)$ is an involutive matrix ($K^2$ is the $2p\times 2p$ identity
matrix) with eigenvalues $1$ and $-1$ of multiplicity $p$ each. The domain
$\cY$ in $\mR^m$ where variable $y$ ranges is assumed to contain the origin
and to be invariant under the linear involution $y\mapsto-y$ (the reflection
with respect to the origin). The neighborhood $\cO_{2p}(0)$ where variable $z$
ranges is supposed to be invariant under the linear involution $z\mapsto Kz$.
The systems we will deal with depend
\begin{itemize}
\item on an angle variable $X\in\mT^N$ ($N\geq 0$) subject to the equation
$\dot{X}=\Omega$ with a fixed Diophantine vector $\Omega\in\mR^N$ (according
to the discussion above, for positive $N$ this actually means a quasi-periodic
dependence on time with frequency vector $\Omega$),
\item on an external parameter $\nu\in\cN\subset\mR^s$ ($s\geq 0$) where $\cN$
is an open domain,
\item and on a small perturbation parameter $\vare\geq 0$.
\end{itemize}
These systems will be assumed to be reversible with respect to the involution
$\cG:(x,y,z,X)\mapsto(-x,-y,Kz,-X)$. We are interested in reducible
quasi-periodic invariant $(n+N)$-tori in such systems in the ``extended phase
space'' $\mT^n\times\cY\times\cO_{2p}(0)\times\mT^N$. It is clear that
\[
P=\dim\Fix\cG=p, \quad Q=\codim\Fix\cG=n+m+p+N, \quad Q-(n+N)=m+p,
\]
so that $P<Q-(n+N)$ since $m\geq 1$, and the situation pertains to the
reversible context~2.

Consider a family of $\cG$-reversible systems on
$\mT^n\times\cY\times\cO_{2p}(0)\times\mT^N$ of the form
\begin{equation}
\begin{aligned}
\dot{x} &= H(y,\nu)+f^\sharp(x,y,z,\nu)+\vare f(x,y,z,X,\nu,\vare), \\
\dot{y} &= \Xi(y,\nu)+g^\sharp(x,y,z,\nu)+\vare g(x,y,z,X,\nu,\vare), \\
\dot{z} &= \Lambda(y,\nu)z+h^\sharp(x,y,z,\nu)+\vare h(x,y,z,X,\nu,\vare), \\
\dot{X} &= \Omega
\end{aligned}
\label{eqsystem}
\end{equation}
(with a $2p\times 2p$ matrix-valued function $\Lambda$), where
$f^\sharp=\Landau{|z|}$, $g^\sharp=\Landau{|z|^2}$, and
$h^\sharp=\Landau{|z|^2}$. Reversibility of~\eqref{eqsystem} with respect to
$\cG$ means that
\begin{gather*}
H(-y,\nu)\equiv H(y,\nu), \quad \Xi(-y,\nu)\equiv\Xi(y,\nu), \\
\Lambda(-y,\nu)K\equiv-K\Lambda(y,\nu)
\end{gather*}
and
\begin{align*}
f^\sharp(-x,-y,Kz,\nu) &\equiv f^\sharp(x,y,z,\nu), &
f(-x,-y,Kz,-X,\nu,\vare) &\equiv f(x,y,z,X,\nu,\vare), \\
g^\sharp(-x,-y,Kz,\nu) &\equiv g^\sharp(x,y,z,\nu), &
g(-x,-y,Kz,-X,\nu,\vare) &\equiv g(x,y,z,X,\nu,\vare), \\
h^\sharp(-x,-y,Kz,\nu) &\equiv -Kh^\sharp(x,y,z,\nu), &
h(-x,-y,Kz,-X,\nu,\vare) &\equiv -Kh(x,y,z,X,\nu,\vare).
\end{align*}
All the functions $H$, $\Xi$, $\Lambda$, $f^\sharp$, $g^\sharp$, $h^\sharp$,
$f$, $g$, and $h$ are assumed to be analytic in all their arguments.

The eigenvalues of the matrix $\Lambda(0,\nu)$ anti-commuting with $K$ come in
pairs $\lambda,-\lambda$ for any $\nu\in\cN$ \cite{S86,S92}. Let the spectrum
of $\Lambda(0,\nu)$ be \emph{simple} for any $\nu\in\cN$ and have the form
\begin{equation}
\begin{gathered}
\pm\alpha_1(\nu),\ldots,\pm\alpha_{d_1}(\nu), \qquad
\pm\rmi\beta_1(\nu),\ldots,\pm\rmi\beta_{d_2}(\nu), \\
\pm\alpha_{d_1+1}(\nu)\pm\rmi\beta_{d_2+1}(\nu),\ldots,
\pm\alpha_{d_1+d_3}(\nu)\pm\rmi\beta_{d_2+d_3}(\nu),
\end{gathered}
\label{eqspec}
\end{equation}
where the numbers $d_1\geq 0$, $d_2\geq 0$, $d_3\geq 0$ do not depend on $\nu$
($d_1+d_2+2d_3=p$), $\alpha_k(\nu)>0$ for all $1\leq k\leq d_1+d_3$,
$\nu\in\cN$, and $\beta_l(\nu)>0$ for all $1\leq l\leq d_2+d_3$, $\nu\in\cN$.

Fix an \emph{arbitrary} (possibly, empty) subset of indices
\[
\fZ\subset\{1;2;\ldots;d_1+d_3\}
\]
consisting of $\kappa$ elements ($0\leq\kappa\leq d_1+d_3$).

\begin{thm}\label{main}
Suppose that
\begin{itemize}
\item $s\geq n+m+d_2+d_3+\kappa$,
\item $\Xi(0,\nu^0)=0$ for some $\nu^0\in\cN$ \kvadra{we will use the notation
$\omega=H(0,\nu^0)\in\mR^n$, $\alpha^0=\alpha(\nu^0)\in\mR^{d_1+d_3}$,
$\beta^0=\beta(\nu^0)\in\mR^{d_2+d_3}$},
\item the vectors $\omega$, $\Omega$, and $\beta^0$ satisfy the following
Diophantine condition\textup{:} there exist constants $\tau>n+N-1$ and
$\gamma>0$ such that
\[
\bigl| \langle j,\omega\rangle+\langle J,\Omega\rangle
+\langle q,\beta^0\rangle \bigr| \geq \gamma\bigl(|j|+|J|\bigr)^{-\tau}
\]
for all $(j,J)\in\mZ^{n+N}\setminus\{0\}$ and $q\in\mZ^{d_2+d_3}$,
$|q|=|q_1|+\cdots+|q_{d_2+d_3}|\leq 2$,
\item the mapping from $\cN$ to $\mR^{n+m+d_2+d_3+\kappa}$ given by
\begin{equation}
\nu\mapsto\bigl( H(0,\nu), \, \Xi(0,\nu), \, \beta(\nu), \,
\alpha_k(\nu), k\in\fZ \bigr)
\label{eqmap}
\end{equation}
is \emph{submersive} at point $\nu^0$, i.e.,
\[
\left. \rank\frac{\partial\bigl( H(0,\nu), \, \Xi(0,\nu), \, \beta(\nu), \,
\alpha_k(\nu), k\in\fZ \bigr)}{\partial\nu} \right|_{\nu=\nu^0} =
n+m+d_2+d_3+\kappa.
\]
\end{itemize}
Then for sufficiently small $\vare$ there exists an
$(s-n-m-d_2-d_3-\kappa)$-dimensional analytic surface $\cS_\vare\subset\cN$
such that for any $\nu\in\cS_\vare$, system~\eqref{eqsystem} admits an
analytic reducible invariant $(n+N)$-torus carrying Diophantine quasi-periodic
motions with frequency vector $(\omega,\Omega)$. The Floquet exponents of this
torus are
\begin{equation}
\begin{gathered}
\underbrace{0,\ldots,0}_m\,, \qquad
\pm\alpha'_1(\nu,\vare),\ldots,\pm\alpha'_{d_1}(\nu,\vare), \qquad
\pm\rmi\beta^0_1,\ldots,\pm\rmi\beta^0_{d_2}, \\
\pm\alpha'_{d_1+1}(\nu,\vare)\pm\rmi\beta^0_{d_2+1},\ldots,
\pm\alpha'_{d_1+d_3}(\nu,\vare)\pm\rmi\beta^0_{d_2+d_3}
\end{gathered}
\label{eqFloq}
\end{equation}
\kvadra{cf.~\eqref{eqspec}}, where $\alpha'_k(\nu,\vare)>0$ for all
$1\leq k\leq d_1+d_3$, $\nu\in\cS_\vare$, and
$\alpha'_k(\nu,\vare)\equiv\alpha^0_k$ for $k\in\fZ$. These tori and the
numbers $\alpha'_k(\nu,\vare)$, $k\notin\fZ$, depend analytically on
$\nu\in\cS_\vare$ and on $\vare^{1/2}$. At $\vare=0$, the surface $\cS_0$
contains $\nu^0$ and all the tori are $\{y=0, \, z=0\}$.
\end{thm}

This theorem describes the persistence of the unperturbed reducible invariant
torus $\{y=0, \, z=0; \, \nu=\nu^0\}$ with the preservation of
\begin{itemize}
\item the frequencies $\omega_1,\ldots,\omega_n$,
\item all the imaginary parts $\pm\beta^0_1,\ldots,\pm\beta^0_{d_2+d_3}$ of
the Floquet exponents,
\item and an arbitrary subcollection (of length $\kappa$) of the pairs of the
real parts $\pm\alpha^0_1,\ldots,\pm\alpha^0_{d_1+d_3}$ of the Floquet
exponents.
\end{itemize}
The situation resembles the so called \emph{partial preservation of Floquet
exponents} \cite{S07Stek} in the ``well developed'' autonomous contexts of KAM
theory (the reversible context~1, the Hamiltonian context, the volume
preserving context, and the dissipative context).

The autonomous case of Theorem~\ref{main} ($N=0$) is the main result of our
previous paper \cite{Snew}. It is tempting to reduce Theorem~\ref{main} with
$N>0$ to this particular case via regarding $(x,X)$ as a new variable $x$ and
$n+N$ as a new dimension $n$. Unfortunately, such an attempt would lead to the
mapping
\begin{equation}
\nu\mapsto\bigl( H(0,\nu), \, \Omega, \, \Xi(0,\nu), \, \beta(\nu), \,
\alpha_k(\nu), k\in\fZ \bigr)
\label{eqbadmap}
\end{equation}
[cf.~\eqref{eqmap}] which cannot be submersive because its second component is
a constant. Moreover, in any case, this approach would require
$s\geq n+N+m+d_2+d_3+\kappa$.

However, one really can reduce Theorem~\ref{main} with arbitrary $N$ to its
particular case $N=0$. Proofs in KAM theory are generally believed to be very
complicated and tedious. Nevertheless, this theory also includes many powerful
methods that enable one to deduce various statements from simpler ones in a
very straightforward manner. Here are some examples (the references we give
just illustrate the methods in question).
\begin{itemize}
\item KAM-type theorems for vector fields and for diffeomorphisms \cite{D82}.
\item A reduction of ``local'' theorems (concerning invariant tori near
equilibria or closed trajectories) to ``global'' ones \cite{S86}.
\item Easy proofs of excitation of elliptic normal modes using
``conventional'' theorems with R\"ussmann-like nondegeneracy conditions
\cite{BHS96LNM,S95Cha}.
\item Herman's method of reducing KAM theorems with weak nondegeneracy
conditions to theorems with nondegeneracy conditions of the submersivity type
\cite{BHS96Gro,BHS96LNM,BS10,S95Cha,S06,S07DCDS,S07Stek}.
\item The easy proofs of the autonomous case of Theorem~\ref{main} (and some
other results in the reversible KAM theory) employing Moser's modifying terms
theorem \cite{M67,S11,Snew}.
\end{itemize}

A reduction of Theorem~\ref{main} with positive $N$ to the autonomous case is
a one more example (probably the simplest one). Introduce an artificial
additional external parameter $\Theta\in\cO_N(0)$ and replace the last
equation $\dot{X}=\Omega$ in systems~\eqref{eqsystem} with the equation
$\dot{X}=\Omega+\Theta$. This does lead to the autonomous framework with
\begin{itemize}
\item $(x,X)$ playing the role of $x$,
\item $n+N$ playing the role of $n$,
\item $(\omega,\Omega)$ playing the role of $\omega$,
\item $(\nu,\Theta)$ playing the role of $\nu$,
\item $s+N$ playing the role of $s$,
\item $\cN\times\cO_N(0)$ playing the role of $\cN$,
\item $(\nu^0,0)\in\mR^{s+N}$ playing the role of $\nu^0$.
\end{itemize}
Indeed, instead of~\eqref{eqmap}, one now has to consider the mapping from
$\cN\times\cO_N(0)$ to $\mR^{n+N+m+d_2+d_3+\kappa}$ given by
\begin{equation}
(\nu,\Theta)\mapsto\bigl( H(0,\nu), \, \Omega+\Theta, \, \Xi(0,\nu), \,
\beta(\nu), \, \alpha_k(\nu), k\in\fZ \bigr)
\label{eqgoodmap}
\end{equation}
[cf.~\eqref{eqbadmap}]. Since the mapping~\eqref{eqmap} is submersive at point
$\nu=\nu^0$, the mapping~\eqref{eqgoodmap} is submersive at point $\nu=\nu^0$,
$\Theta=0$.

Apply the autonomous Theorem~\ref{main} proven in \cite{Snew} using Moser's
modifying terms theorem \cite{M67}. We arrive at the conclusion that for
sufficiently small $\vare$, there exists an
$(s-n-m-d_2-d_3-\kappa)$-dimensional analytic surface
$\cS_\vare\subset\bigl(\cN\times\cO_N(0)\bigr)$ such that for any
$(\nu,\Theta)\in\cS_\vare$, the modified system~\eqref{eqsystem} [with
$\dot{X}=\Omega+\Theta$ instead of $\dot{X}=\Omega$] admits an analytic
reducible invariant $(n+N)$-torus carrying Diophantine quasi-periodic motions
with frequency vector $(\omega,\Omega)$. The Floquet exponents of this torus
are
\begin{equation}
\begin{gathered}
\underbrace{0,\ldots,0}_m\,, \qquad
\pm\alpha'_1(\nu,\Theta,\vare),\ldots,\pm\alpha'_{d_1}(\nu,\Theta,\vare),
\qquad \pm\rmi\beta^0_1,\ldots,\pm\rmi\beta^0_{d_2}, \\
\pm\alpha'_{d_1+1}(\nu,\Theta,\vare)\pm\rmi\beta^0_{d_2+1},\ldots,
\pm\alpha'_{d_1+d_3}(\nu,\Theta,\vare)\pm\rmi\beta^0_{d_2+d_3}
\end{gathered}
\label{eqnew}
\end{equation}
[cf.~\eqref{eqFloq}], where $\alpha'_k(\nu,\Theta,\vare)>0$ for all
$1\leq k\leq d_1+d_3$, $(\nu,\Theta)\in\cS_\vare$, and
$\alpha'_k(\nu,\Theta,\vare)\equiv\alpha^0_k$ for $k\in\fZ$. These tori and
the numbers $\alpha'_k(\nu,\Theta,\vare)$, $k\notin\fZ$, depend analytically
on $(\nu,\Theta)\in\cS_\vare$ and on $\vare^{1/2}$. At $\vare=0$, the surface
$\cS_0$ contains the point $(\nu^0,0)$ and all the tori are $\{y=0, \, z=0\}$.

Now it suffices to verify that the surface $\cS_\vare$ lies in fact in
$\cN\times\{0\}$ for each $\vare$, i.e.\ that $\Theta=0$ on $\cS_\vare$ [and
the dependence of $\alpha'$ on $\Theta$ in~\eqref{eqnew} is dummy]. Consider
the $(n+N)$-torus $\cT$ corresponding to an arbitrary point
$(\nu,\Theta)\in\cS_\vare$. According to the standard reflection lemma of
Section~\ref{revers}, one can introduce in $\cT$ a coordinate frame
$\varp\in\mT^n$, $\Phi\in\mT^N$ in which the dynamics on $\cT$ takes the form
\begin{equation}
\dot{\varp}=\omega, \quad \dot{\Phi}=\Omega
\label{eqoO}
\end{equation}
and the restriction of the involution $\cG:(x,X,y,z)\mapsto(-x,-X,-y,Kz)$ to
$\cT$ takes the form
\begin{equation}
\cG|_{\cT}:(\varp,\Phi)\mapsto(-\varp,-\Phi).
\label{eqcG}
\end{equation}
The torus $\cT$ is close to the unperturbed $(n+N)$-torus $\{y=0, \, z=0\}$ at
$\nu=\nu^0$, $\Theta=0$, $\vare=0$ for which one can set $\varp=x$, $\Phi=X$.
Consequently, the coordinates $\varp$, $\Phi$ in $\cT$ can be chosen in such a
way that the torus $\cT$ will be given by the equations
\[
x=\varp+A(\varp,\Phi), \quad X=\Phi+B(\varp,\Phi), \quad
y=C(\varp,\Phi), \quad z=D(\varp,\Phi),
\]
$A$, $B$, $C$, and $D$ being small analytic functions with values in $\mR^n$,
$\mR^N$, $\mR^m$, and $\mR^{2p}$, respectively. In view of~\eqref{eqcG} these
functions satisfy the identities
\begin{align*}
A(-\varp,-\Phi) &\equiv -A(\varp,\Phi),
& B(-\varp,-\Phi) &\equiv -B(\varp,\Phi), \\
C(-\varp,-\Phi) &\equiv -C(\varp,\Phi),
& D(-\varp,-\Phi) &\equiv KD(\varp,\Phi).
\end{align*}
Differentiating the relation $X=\Phi+B(\varp,\Phi)$ with respect to time and
taking~\eqref{eqoO} and the equation $\dot{X}=\Omega+\Theta$ into account, one
obtains
\[
\Omega+\Theta\equiv\Omega+\frac{\partial B}{\partial\varp}\omega+
\frac{\partial B}{\partial\Phi}\Omega,
\]
i.e.,
\begin{equation}
\frac{\partial B}{\partial\varp}\omega+\frac{\partial B}{\partial\Phi}\Omega
\equiv\Theta.
\label{eqyes}
\end{equation}
The mean value of the left-hand side of~\eqref{eqyes} over $\mT^{n+N}$
vanishes, whence $\Theta=0$. Moreover, the identity~\eqref{eqyes} with
$\Theta=0$ and the fact that the numbers
$\omega_1,\ldots,\omega_n,\Omega_1,\ldots\Omega_N$ are independent over
rationals imply easily that $B=\const$ (cf.~\cite[section~6.2]{S07DCDS}).
Since the function $B$ is odd, $B\equiv 0$. The proof of Theorem~\ref{main}
with arbitrary $N\geq 0$ is completed.

Remark finally that the reversible KAM theory can be extended to systems
quasi-periodic not only in time but also in some spatial variables
\cite{HS10}.

\section*{Acknowledgments}

The author is indebted to H.W.~Broer for a long-term friendship and
collaboration and many fruitful discussions on quasi-periodicity in dynamical
systems. This study was partially supported by a grant of the President of the
Russia Federation, project No.\ NSh-8462.2010.1.

\end{document}